\documentclass[12pt]{article}
\usepackage{amsfonts}
\usepackage{hyperref}
\usepackage{mathrsfs}
\usepackage{indentfirst}
\usepackage{amsmath}
\usepackage{amssymb}
\usepackage[amsmath,thmmarks]{ntheorem}
\usepackage{graphicx}
\usepackage{tikz}

\newtheorem{definition}{\bf Definition}[section]

\newtheorem{theorem}{\bf Theorem}[section]
\newtheorem{remark}{\bf Remark}[section]
\newtheorem{corollary}{\bf Corollary}[section]
\newtheorem{example}{\bf Example}[section]
\newtheorem{algorithm}{\bf Algorithm}[section]

\newtheorem{proposition}{\bf Proposition}[section]

\begin{document}
\setcounter{page}{1}

\title{{\textbf{Minimax programming problems subject to addition-{\L}ukasiewicz fuzzy relational inequalities and their optimal solutions}\thanks {Supported by
the National Natural Science Foundation of China (No.12071325)}}}
\author{Xue-ping Wang\footnote{Corresponding author. xpwang1@hotmail.com; fax: +86-28-84761502}, Meng Li\footnote{\emph{E-mail address}: 1228205272@qq.com}, Qian-yu Shu\footnote{\emph{E-mail address}: 34956229@qq.com}\\
\emph{School of Mathematical Sciences, Sichuan Normal University,}\\
\emph{Chengdu 610066, Sichuan, People's Republic of China}}

\newcommand{\pp}[2]{\frac{\partial #1}{\partial #2}}
\date{}
\maketitle

\begin{quote}
{\bf Abstract} This article focuses on minimax programming problems subject to addition-{\L}ukasiewicz fuzzy relational inequalities. We first establish two necessary and sufficient conditions that a solution of the fuzzy relational inequalities is a minimal one and explore the existence condition of the unique minimal solution. We also supply an algorithm to search for minimal solutions of the fuzzy relational inequalities starting from a given solution. We then apply minimal solutions of the fuzzy relational inequalities to the minimax programming problems for searching optimal solutions. We provide two algorithms to solve a kind of single variable optimization problems, and obtain the greatest optimal solution. The algorithm for finding minimal solutions of a given solution are also used for searching minimal optimal solutions.

{\textbf{\emph{Keywords}}:}\, Fuzzy relational inequality; Addition-{\L}ukasiewicz composition; Minimal solution; Minimax programming problem; Optimal solution.\\
\end{quote}

\section{Introduction}
A fuzzy relation plays a central role in fuzzy system, and a fuzzy relational equation is one of the main topics of fuzzy relations. In 1976, Sanchez \cite{Sanchez76} first introduced and investigated fuzzy relational equations with max-min composition. Since then, many researchers all over the world have discussed a variety of fuzzy relational equations (inequalities) with different composition operators, such as max-min, max-product, max-{\L}ukasiewicz, addition-min, etc., see e.g. \cite{DiNola89,Li2012,Martin2017,Peeva04}. They have tried to explore an effective method to describe the solution set of fuzzy relational equations (inequalities) by their minimal solutions, unique solution, approximate solutions and so on, respectively, see e.g. \cite{Ab06,DiNola89,Li2012,Peeva04,Qiu}.
It is well known that the data transmission mechanism in BitTorrent-like peer-to-peer (P2P) file sharing system has been reduced to a system of fuzzy relational inequalities with addition-min composition as below.
\begin{equation}\label{eq1}
\left\{\begin{array}{ll}
a_{11}\wedge x_1 + a_{12}\wedge x_{2} + \cdots + a_{1n}\wedge x_n\geq b_1,\\
a_{21}\wedge x_1 + a_{22}\wedge x_{2} + \cdots + a_{2n}\wedge x_n\geq b_2,\\
\cdots\\
a_{m1}\wedge x_1 + a_{m2}\wedge x_{2} + \cdots + a_{mn}\wedge x_n\geq b_m,
\end{array}\right.
\end{equation}
 where $a_{ij},x_{j} \in [0,1]$, $b_i > 0$, $a_{ij} \wedge x_{j} = \min \{a_{ij}, x_j\}$ with $i \in \{1, 2, \cdots, m\}$ and $j \in \{1, 2, \cdots, n\}$, and the operation $``+"$ is the ordinary addition \cite{Li2012,Yang2014}. In such a system, there are $n$ users $ A_1, A_2, \cdots , A_n$ who simultaneously download some file data. The $jth$ user $A_j$ delivers the file data with quality level $x_j$ to $A_i$, and $a_{ij}$ is the communicative bandwidth between $A_i$ and $A_j$. So that $a_{ij}\wedge x_j$ is the actual network traffic that $A_i$ receives from $A_j$, and the communicative quality requirement of download traffic of $A_i$ is at least $b_i$ with $b_i > 0$, then $A_i$ receiving the file data from other users is $$a_{i1}\wedge x_1 + \cdots +a_{i{}i-1}\wedge x_{i-1}+a_{i{i+1}}\wedge x_{i+1}+\cdots+ a_{in}\wedge x_n\geq b_i. $$ Based on the above, we add the $a_{ii}\wedge x_i$ in the $i$th inequality and get the system (\ref{eq1}).

  Li and Yang \cite{Li2012} gave a sufficient and necessary condition for the system (\ref{eq1}) to be solvable, and supplied a method to get one minimal solution. Yang \cite{Yang2018} showed that the solution set of system (\ref{eq1}) is convex and has an infinite number of minimal solutions, unlike the solution sets of most of the fuzzy relational inequalities which have a finite number of minimal solutions. To avoid a network congestion and ensure the data transmission, many scholars have studied optimization models of different objective functions with the system (\ref{eq1}) as constraints, such as linear programming \cite{Guo,Guu2,Yang2014}, multilevel programming \cite{Guu19,Yang2015}, and min-max programming \cite{chi,Lin20,Wu1,Wu2,Wu3,Yang2017,Yang2016}. In particular, Yang \cite{Yang2014} proved that all the optimal solutions of the related optimization models are minimal solutions of system (\ref{eq1}). This inspires that many authors do their best to discuss the existence of minimal solutions and try to find minimal solutions of system (\ref{eq1}), see e.g. \cite{Li2012,Li2021,Wu,Yang2014,SJYang,Yang2018}. Fortunately, Li and Wang \cite{Li2022} presented an algorithm to compute all minimal solutions of system (\ref{eq1}).

One could imagine that the subject of solving the system (\ref{eq1}) has already been dealt with, however, this is not the case. Indeed, the operator $\wedge$ in the system (\ref{eq1}) is just one of the continuous triangular norms (t-norms for short), and we usually denote $\wedge$ by $T_M$. From the theory of t-norms, we know that there are three elementary kinds of typical continuous t-norms, i.e., minimum ($T_M$), product ($T_P$) and {\L}ukasiewicz t-norm ($T_L$), and starting with the three t-norms, it is possible to construct all continuous t-norms by means of isomorphic transformations and ordinal sums \cite{LING1965}. Therefore, from the continuous t-norm point of view, we just need to consider three kinds of mathematical models for the P2P file sharing systems, i.e., the system with addition-min composition, the system with addition-product composition and the system with addition-{\L}ukasiewicz composition. For the first system, i.e., the system (\ref{eq1}), Li and Wang \cite{Li2022} supplied an algorithm to calculate all its minimal solutions, therefore, we can describe the solution set of system (\ref{eq1}) by all minimal solutions and its greatest solution when it is solvable. Moreover, if we replace $\wedge$ in the system (\ref{eq1}) by $T_P$ then the system (\ref{eq1}) can be transformed into a classically linear inequalities assigned on [0,1] whose solution set can be characterized by a current method \cite{So}. Therefore, an interesting question is how to describe the solution set of system (\ref{eq1}) when we replace $\wedge$ in the system (\ref{eq1}) by $T_L$. This motivates us for studying minimal solutions of the following mathematical model for the P2P file sharing system:
 \begin{equation}\label{eq2}
\left\{\begin{array}{ll}
T_L(a_{11}, x_1) + T_L(a_{12}, x_{2}) + \cdots + T_L(a_{1n}, x_n)\geq b_1,\\
T_L(a_{21}, x_1) +T_L(a_{22}, x_{2}) + \cdots + T_L(a_{2n}, x_n)\geq b_2,\\
\cdots\\
T_L(a_{m1}, x_1) +T_L(a_{m2}, x_{2}) + \cdots + T_L(a_{mn}, x_n)\geq b_m,
\end{array}\right.
\end{equation}
where $a_{ij},x_{j} \in [0,1]$, $b_i > 0$, $T_L(a_{ij}, x_{j}) = \max \{a_{ij}+x_j-1, 0\}$ with $i \in \{1, 2, \cdots, m\}$ and $j \in \{1, 2, \cdots, n\}$, and the operation $``+"$ is the ordinary addition. To avoid the network congestion and improve the stability of data transmission, we further apply their minimal solutions of system \eqref{eq2} to search the optimal solutions of the following minimax programming problem:
\begin{align}\label{eq6}
\mbox{Minimize}\quad&Z(x)= x_1\vee x_2\vee \cdots\vee x_n \nonumber\\
\mbox{subject to}\quad& \left\{\begin{array}{ll}
T_L(a_{11}, x_1) + T_L(a_{12}, x_2) + \cdots + T_L(a_{1n}, x_n)\geq b_1,\\
T_L(a_{21}, x_1) +T_L(a_{22}, x_2) + \cdots + T_L(a_{2n}, x_n)\geq b_2,\\
\cdots\\
T_L(a_{m1}, x_1) +T_L(a_{m2}, x_2) + \cdots + T_L(a_{mn}, x_n)\geq b_m.
\end{array}\right.
\end{align}

The remainder of this article is organized by five sections. Section 2 presents some definitions and properties of system (\ref{eq2}). Section 3 supplies a necessary and sufficient condition for a solution of system (\ref{eq2}) to be a minimal one, and investigates the existence condition of the unique minimal solution. In Section 4, we first show that for any solution of system (\ref{eq2}) there is a minimal one that is less than or equal to the solution, and then propose an algorithm to search for minimal solutions of system (\ref{eq2}) starting from a fixed solution with computational complexity $O(m^3)$ or $O(n^3)$. In Section 5, we first transform the problem \eqref{eq6} into a single variable optimization problem, and then provide two algorithms for solving the single variable optimization problem and obtain the greatest optimal solution of problem \eqref{eq6}. We further apply the algorithm established in Section 4 to find minimal optimal solutions of problem \eqref{eq6}. A conclusion is given in Section 6.

\section{Preliminaries}
In this section, we give basic definitions and properties of system (\ref{eq2}).
\begin{definition}[\cite{Kl}]\label{de1.1}
\emph{A binary operator $T: [0,1]^2 \rightarrow [0,1]$ is said to be a t-norm,
if it is commutative, associative, increasing with respect to each variable and has an identity $1$, i.e.,  $T(x, 1) = x$ for all $x\in [0,1]$.}
\end{definition}

Let $I = \{1, 2, \cdots , m\}$ and $J = \{ 1, 2, \cdots, n \}$ be two index sets. For $x^1 = (x_1^1, x_2^1, \cdots, x_n^1)$, $x^2 = (x_1^2, x_2^2, \cdots, x_n^2) \in [0,1]^{n}$, define ${x^1} \le {x^2}$ iff ${x_j}^1 \le {x_j}^2$ for arbitrary $j \in J$; and define ${x^1} < {x^2}$ iff ${x_j}^1 \le {x_j}^2$ and there is a ${j_0} \in J$ such that $x_{j_0}^1 < x_{j_0}^2$. We further define $x^1\vee x^2$ by  $x^1\vee x^2=(x_1^1\vee x_1^2, x_2^1\vee x_2^2, \cdots, x_n^1\vee x_n^2)$. Then the system (\ref{eq2}) can be tersely described as follows:
$$\sum_{j\in J}T_L(a_{ij}, x_{j})\geq b_i \mbox{ for all }i\in I$$
or
$$A \odot {x^T} \geq {b^T}$$
where $A=(a_{ij})_{m \times n}$, $x=(x_1, x_2, \cdots, x_n)$, $b=(b_1, b_2, \cdots, b_m)$ and
$(a_{i1}, a_{i2}, \cdots , a_{in}) \odot ( x_1, x_2, \cdots, x_n)^T =T_L(a_{i1}, x_1) + T_L(a_{i2}, x_{2}) + \cdots + T_L(a_{in}, x_n)\geq b_i$. The system (\ref{eq2}) is called solvable if there exists an $x \in [0,1]^{n}$ satisfying $A \odot {x^T} \geq {b^T}$. Denote the set of all the solutions of system (\ref{eq2}) by $$S(A,b) = \{x \in [0,1]^{n} | A \odot x^T \geq b^T\},$$ and $U\setminus V=\{x\in U \mid x\notin V\}$ where $U$ and $V$ are two sets.

\begin{definition}[\cite{DiNola89}]
 \emph{An $\hat x \in S(A,b)$ is said to be the greatest solution if $x \le \hat x$ for all $x \in S(A,b)$; an $\check{x}\in S(A,b)$ is said to be a minimal solution if $x \le \check{x}$ for any $x \in S(A,b)$, implies $x = \check{x}$}.
\end{definition}

From Definition \ref{de1.1}, the following is clear.
\begin{theorem}\label{th2.1}
$S(A,b)\neq \emptyset$ iff $\sum_{j\in J}a_{ij}\geq b_i$ for all $i\in I$.
\end{theorem}

Theorem \ref{th2.1} and Definition \ref{de1.1} imply the following proposition.
\begin{proposition}\label{pro2.1}
For the system (\ref{eq2}), we have:
\item [(i)] If $S(A, b)\neq \emptyset$, then $(1,1,\cdots, 1)$ is its greatest solution.
\item [(ii)] Let $x\in S(A,b)$, $x^*\in [0, 1]^n$. Then $x\le x^*$ implies $x^*\in S(A,b)$.
\end{proposition}

 According to Proposition \ref{pro2.1} (ii), the following proposition holds.
\begin{proposition}\label{pro2.2}
If $x^1, x^2\in S(A,b)$, then $x^1\vee x^2\in S(A,b)$.
\end{proposition}
\begin{theorem}\label{th2.2}
If $S(A,b)\neq \emptyset$, and there is an $i\in I$ such that $\sum_{j\in J}a_{ij}= b_i$ and $a_{ij}\neq 0$ for any $j\in J$, then $(1,1,\cdots, 1)$ is the unique solution of system (\ref{eq2}).
\end{theorem}
\begin{proof}
Since $S(A,b)\neq \emptyset$, we have $(1,1,\cdots, 1)\in S(A,b)$. Let $x = (x_1, x_2, \cdots, x_n)\in S(A,b)$ and $x\neq (1,1,\cdots, 1)$. Then there is a $j\in J$ such that $x_j< 1$. If there is an $i\in I$ such that $\sum_{j\in J}a_{ij}= b_i$ and $a_{ij}\neq 0$ for any $j\in J$, then $ \max\{a_{ij}+x_j-1,0\}<a_{ij}$. Therefore, $\sum_{k\in J\setminus\{j\}}T_L(a_{ik},x_k)+T_L(a_{ij},x_j)<\sum_{j\in J}a_{ij}=b_i$, contrary to $x\in S(A,b)$. This follows that $(1,1,\cdots, 1)$ is the unique solution of system (\ref{eq2}).
\end{proof}
\begin{proposition}
 If $x, y \in S(A, b)$ with $x\geq y$. Then $\lambda x+(1-\lambda)y\in S(A, b)$ for any $\lambda\in [0, 1]$.
\end{proposition}
\begin{proof}
  Since $x, y \in S(A, b)$ with $x\geq y$, for any $\lambda\in [0, 1]$ and $i\in I$, we have
  \begin{eqnarray*} &&\sum_{j\in J}T_L(a_{ij},\lambda x_j+(1-\lambda)y_j)
 \\&=&\sum_{j\in J}\max\{a_{ij}+\lambda x_j+(1-\lambda)y_j-1, 0\}\\&=&\sum_{j\in J}\max\{a_{ij}+y_j-1+\lambda (x_j-y_j), 0\}\\&\geq&\sum_{j\in J}\max\{a_{ij}+y_j-1, 0\}\\&\geq&b_i.
 \end{eqnarray*}
 Therefore, $\lambda x+(1-\lambda)y\in S(A, b)$.
\end{proof}
\section{Minimal solutions of system (\ref{eq2})}
This section shows a sufficient and necessary condition for a solution of system \eqref{eq2} being a minimal one. It then characterizes the uniqueness of minimal solutions to the system (\ref{eq2}).

Let $x = (x_1, x_2, \cdots, x_n)\in [0, 1]^n$. Denote $J^*(x)=\{j\in J|x_j>0\}$.
\begin{theorem}\label{th3.1}
Let $x = (x_1, x_2, \cdots, x_n)\in S(A,b)$. Then $x$ is a minimal solution of system (\ref{eq2}) iff there exists an $i\in I$ such that $\sum_{j\in J}T_L(a_{ij}, x_{j})= b_i$, and $a_{ij}+x_j-1>0$ for any $j\in J^*(x)$.
\end{theorem}
 \begin{proof} Let $x = (x_1, x_2, \cdots, x_n)\in S(A,b)$. Assume that for any $i\in I$ either $\sum_{j\in J}T_L(a_{ij}, x_{j})>b_i$, or there is a $j_0\in J^*(x)$ such that $a_{ij_0}+x_{j_0}-1\leq 0$.

If there is a $j_0\in J^*(x)$ such that $a_{ij_0}+x_{j_0}-1\leq 0$ for any $i\in I$, then $a_{ij_0}+u-1\leq 0$ for all $i\in I$, $u\in [0, x_{j_0}]$. So that if let $u_0\in [0, x_{j_0})$ and define $x' = (x'_1, x'_2, \cdots, x'_n)$ with $$x'_j = \left\{ \begin{array}{l}
u_0, j = j_0,\\
x_j, j \ne j_0.
\end{array} \right.$$
Then it is clear that $x'< x$ and $T_L(a_{ij_0},x_{j_0})=\max\{a_{ij_0}+x_{j_0}-1, 0\}=0=\max\{a_{ij_0}+x'_{j_0}-1, 0\}=T_L(a_{ij_0},x'_{j_0})$ for all $i\in I$.
Hence, for any $i\in I$,
 \begin{eqnarray*} \sum_{j\in J}T_L(a_{ij}, x'_{j})&=&T_L(a_{ij_0},x'_{j_0})+\sum_{j\in J\setminus\{j_0\}}T_L(a_{ij}, x'_{j})\\&=&T_L(a_{ij_0},x_{j_0})+ \sum_{j\in J\setminus\{j_0\}}T_L(a_{ij}, x_{j})\\&=&\sum_{j\in J}T_L(a_{ij}, x_{j})\\&\geq& b_i.
 \end{eqnarray*}
This indicates that $x'\in S(A,b)$, which contradicts that $x$ is a minimal solution of system (\ref{eq2}).

In the case of $\sum_{j\in J}T_L(a_{ij}, x_{j})> b_i$ for all $i\in I$, let $I_1=\{i\in I|a_{ij_0}+x_{j_0}-1>0\}$ with $j_0\in J^*(x)$. Then $I_1\neq \emptyset$. From $\sum_{j\in J}T_L(a_{ij}, x_{j})=T_L(a_{ij_0},x_{j_0})+\sum_{j\in J\setminus\{j_0\}}T_L(a_{ij}, x_{j})=\max\{a_{ij_0}+x_{j_0}-1, 0\}+ \sum_{j\in J\setminus\{j_0\}}T_L(a_{ij}, x_{j})> b_i$ for any $i\in I_1$, we have $a_{ij_0}+x_{j_0}-1> b_i-\sum_{j\in J\setminus\{j_0\}}T_L(a_{ij}, x_{j})$. Hence, $x_{j_0}>1-a_{ij_0}+b_i-\sum_{j\in J\setminus\{j_0\}}T_L(a_{ij}, x_{j})$ for all $i\in I_1$.
Denote $\dot{x}_{j_0}=\max\limits_{i\in I_1}\{0, 1-a_{ij_0}+b_i-\sum_{j\in J\setminus\{j_0\}}T_L(a_{ij}, x_{j})\}$. Define $y = (y_1, y_2, \cdots, y_n)$ with $$y_j = \left\{ \begin{array}{l}
\dot{x}_{j_0}, j = j_0,\\
x_j, j \ne j_0.
\end{array} \right.$$
Then it is clear that $y<x$. We distinguish two cases:\\
Case 1. For any $i\in I_1$,
 \begin{eqnarray*} &&\sum_{j\in J}T_L(a_{ij}, y_{j})\\&=&T_L(a_{ij_0},y_{j_0})+ \sum_{j\in J\setminus\{j_0\}}T_L(a_{ij}, y_{j})
 \\&=&\max\{a_{ij_0}+y_{j_0}-1, 0\}+ \sum_{j\in J\setminus\{j_0\}}T_L(a_{ij}, x_{j})\\&=&\max\{a_{ij_0}+\max\limits_{i\in I_1}\{0,1-a_{ij_0}+b_i-\sum_{j\in J\setminus\{j_0\}}T_L(a_{ij}, x_{j})\}-1, 0\}+ \sum_{j\in J\setminus\{j_0\}}T_L(a_{ij}, x_{j})\\&\geq&\max\{b_i-\sum_{j\in J\setminus\{j_0\}}T_L(a_{ij}, x_{j}),0\}+ \sum_{j\in J\setminus\{j_0\}}T_L(a_{ij}, x_{j})\\&\geq&b_i-\sum_{j\in J\setminus\{j_0\}}T_L(a_{ij}, x_{j})+\sum_{j\in  J\setminus\{j_0\}}T_L(a_{ij},x_j)\\&=&b_i.
 \end{eqnarray*}
Case 2. For any $i\in I\setminus I_1$, since $y_{j_0}< x_{j_0}$, we have $0\leq T_L(a_{ij_0}, y_{j_0})=\max\{a_{ij_0}+y_{j_0}-1, 0\}\leq \max\{a_{ij_0}+x_{j_0}-1, 0\}= T_L(a_{ij_0},x_{j_0})= 0$, i,e., $T_L(a_{ij_0},y_{j_0})=T_L(a_{ij_0},x_{j_0})$. Therefore, $\sum_{j\in J}T_L(a_{ij}, y_{j})=\sum_{j\in J}T_L(a_{ij}, x_{j})\geq b_i$ for all $i\in I\setminus I_1$.

Cases 1 and 2 mean that $y\in S(A,b)$, contrary to the fact that $x$ is a minimal solution of system (\ref{eq2}).

Conversely, suppose that there exists an $i\in I$ such that $\sum_{j\in J}T_L(a_{ij}, x_{j})= b_i$, and $a_{ij}+x_j-1>0$ for any $j\in J^*(x)$. Now, let $z = (z_1, z_2, \cdots, z_n)\in S(A,b)$ with $z< x$. Then there exists a $j\in J$ such that $0\leq z_j<x_j$. It is easily seen that $j\in J^*(x)$, and $$T_L(a_{ij}, z_j)=\max\{a_{ij}+z_j-1, 0\}< \max\{a_{ij}+x_j-1, 0\}=T_L(a_{ij}, x_j).$$
Moreover, $z_k\leq x_k$ for any $k\in J \setminus\{j\}$, which leads to $$T_L(a_{ik}, z_k)=\max\{a_{ik}+z_k-1, 0\}\leq\max\{a_{ik}+x_k-1, 0\}=T_L(a_{ik}, x_k).$$
Therefore, $\sum_{j\in J}T_L(a_{ij}, z_{j})< \sum_{j\in J}T_L(a_{ij}, x_{j})=b_i$, contrary to $z\in S(A,b)$.
\end{proof}

Denote the set of all minimal solutions of system (\ref{eq2}) by $\check{S}(A,b)$.
\begin{proposition}\label{pro3.1}
Let $S(A,b)\neq \emptyset$ and define \begin{eqnarray}\label{eq03}y = (y_1, y_2, \cdots, y_n)\mbox{ with }y_k = \left\{ \begin{array}{l}
\max\limits_{i\in I}\{0, 1+b_i-\sum_{t\in J}a_{it}\}, k = j,\\
1, k \ne j.
\end{array} \right.\end{eqnarray}
Then $y\in S(A,b)$. Furthermore, if $a_{it}>0$ for all $i\in I$, $t\in J$ and $a_{ij}+y_j-1>0$ for any $i\in I$, then $y\in \check{S}(A,b)$.
\end{proposition}
 \begin{proof}
 We first verify that $y\in S(A,b)$. By the construction of $y$, $y_j\geq 1+b_i-\sum_{t\in J}a_{it}$, i.e., $a_{ij}+y_j-1\geq b_i-\sum_{t\in J\setminus\{j\}}a_{it}$ for any $i\in I$.
 It follows that $\max\{a_{ij}+y_j-1,0\}\geq b_i-\sum_{t\in J\setminus\{j\}}a_{it}$. Therefore, $T_L(a_{ij},y_{j})+ \sum_{t\in J\setminus\{j\}}a_{it}\geq b_i$ for any $i\in I$, which indicates $y\in S(A,b)$.

Next, we show that $y\in \check{S}(A,b)$ if $a_{it}>0$ for all $i\in I$, $t\in J$ and $a_{ij}+y_j-1>0$ for any $i\in I$. First note that $y_j>0$. Then there is an $i_0\in I$ such that $y_j= 1+b_{i_0}-\sum_{t\in J}a_{{i_0}t}>0$ and $a_{{i_0}j}+y_j-1>0$. It follows from $T_L(a_{{i_0}t},1)=a_{{i_0}t}$ that
 \begin{eqnarray*}\sum_{t\in J}T_L(a_{{i_0}t}, y_{t})&=& \sum_{t\in J\setminus\{j\}}T_L(a_{{i_0}t}, y_{t})+T_L(a_{i_0j}, y_{j})\\&=&\sum_{t\in J\setminus\{j\}}T_L(a_{{i_0}t}, y_{t})+\max\{a_{{i_0}j}+y_j-1, 0\}\\&=&\sum_{t\in J\setminus\{j\}}a_{{i_0}t}+ a_{i_0j}+y_j-1\\&=&\sum_{t\in J\setminus\{j\}}a_{{i_0}t}+ a_{i_0j}+1+b_{i_0}-\sum_{t\in J}a_{{i_0}t}-1\\&=&b_{i_0}. \end{eqnarray*}

Now, suppose that there is a $\bar{y}= (\bar{y}_1, \bar{y}_2, \cdots, \bar{y}_n)\in S(A,b)$ with $\bar{y}< y$. Then there is a $t_0\in J$ such that $\bar{y}_{t_0}< y_{t_0}$. Since $a_{it}>0$ for all $i\in I$, $t\in J$ and by the construction of $y$, we have $a_{{i_0}{t_0}}+y_{t_0}-1>0$, then $T_L(a_{{i_0}{t_0}}, y_{t_0})=\max\{a_{{i_0}{t_0}}+y_{t_0}-1, 0\}> \max\{a_{{i_0}{t_0}}+\bar{y}_{t_0}-1, 0\}=T_L(a_{{i_0}{t_0}}, \bar{y}_{t_0}).$ Thus
 \begin{eqnarray*}
 b_{i_0} &=&\sum_{t\in J}T_L(a_{{i_0}t}, y_{t})\\&=& \sum_{t\in J\setminus\{t_0\}}T_L(a_{{i_0}t}, y_{t})+T_L(a_{{i_0}{t_0}}, y_{t_0})\\&>&\sum_{t\in  J\setminus\{t_0\}}T_L(a_{{i_0}t}, \bar{y}_{t})+T_L(a_{{i_0}{t_0}},\bar{y}_{t_0}) \\&=&\sum_{t\in J}T_L(a_{{i_0}t}, \bar{y}_t),
 \end{eqnarray*}
contrary to $\bar{y}\in S(A,b)$. Therefore, $y\in \check{S}(A,b)$.
 \end{proof}

The following example illustrates Proposition \ref{pro3.1}.
\begin{example}\label{ex1}
\emph{Consider the following fuzzy relational inequalities:}
\begin{equation*}
\left\{\begin{array}{ll}
T_L(0.5, x_1) + T_L(0.9, x_{2}) + T_L(0.7, x_3)\geq 1.7,\\
T_L(0.7, x_1) +T_L(0.5, x_{2}) + T_L(0.6, x_3)\geq 1.2,\\
T_L(0.6, x_1) +T_L(0.8, x_{2})+ T_L(0.9, x_3)\geq 1.8.
\end{array}\right.
\end{equation*}
\end{example}
Obviously, $(1, 1, 1)\in S(A,b)$. A simple calculation leads to  \begin{eqnarray*}y_1&=&\max\limits_{i\in \{1,2,3\}}\{0, 1+b_i-\sum_{j\in \{1,2,3\}}a_{ij}\}\\&=&\max\{0.6, 0.4, 0.5\}\\&=&0.6,\end{eqnarray*} and $a_{i1}+0.6-1>0$ for all $i\in \{1,2,3\}$. Then by Proposition \ref{pro3.1}, $y=(0.6, 1, 1)$ is a minimal solution.

Similarly, both $(1, 0.6, 1)$ and $(1, 1, 0.6)$ are minimal solutions.

From Example \ref{ex1}, minimal solutions of system (\ref{eq2}) are usually not unique. So that a quite natural problem is: what is the condition for a minimal solution of system (\ref{eq2}) being unique? Next, we will explore the condition under which the system (\ref{eq2}) has a unique minimal solution.

Let $x = (x_1, x_2, \cdots, x_n)\in S(A,b)$, and for any $j\in J$, define \begin{equation}\label{eq3} F_{j}(x) = \{x_{*j} \in [0,1] |T_L(a_{ij}, x_{*j}) + \sum\limits_{k \in J\setminus\{j\}}T_L( a_{ik}, x_k) \ge b_i \mbox{ for any } i \in I\}\end{equation} and \begin{equation}\label{eq4} \delta_j(x)=\min\{x_{*j}|x_{*j}\in F_{j}(x)\}.\end{equation}

From (\ref{eq3}), we deduce the following proposition.
\begin{proposition}\label{pro3.2}
Let $x,y\in S(A,b)$ with $x \le y$. Then $F_j(x)\subseteq F_j(y)$ for any $j\in J$.
\end{proposition}

In what follows, let $\hat x=(1,1,\cdots, 1)$. Then we have the following theorem.
\begin{theorem}\label{th3.2}
If $S(A,b) \ne \emptyset$, then $\{x_{j}|x = (x_1, x_2, \cdots, x_n)\in S(A,b)\} = F_{j}(\hat{x})$.
\end{theorem}
\begin{proof}
If $S(A,b) \ne \emptyset$, then $x\le \hat{x}$ for any $x\in S(A,b)$. According to Proposition \ref{pro3.2}, it is obvious that $\{x_{j}|x = (x_1, x_2, \cdots, x_n)\in S(A,b)\} \subseteq F_{j}(\hat{x})$.

 Now, let $x_{j} \in F_{j}(\hat{x})$ and define $y = (y_1, y_2, \cdots, y_n)$ with $$y_k = \left\{ \begin{array}{l}
1, k\ne j,\\
x_{j}, k = j.
\end{array} \right.$$
Easily see that $y\in S(A,b)$, then $x_{j} \in \{x_{j}|x = (x_1, x_2, \cdots, x_n)\in S(A,b)\}$, i.e., $F_{j}(\hat{x}) \subseteq \{x_{j}|x = (x_1, x_2, \cdots, x_n)\in S(A,b)\}$. Therefore, $\{x_{j}|x = (x_j)_{j \in J} \in S(A,b)\} = F_{j}(\hat{x})$.
\end{proof}
\begin{corollary}\label{cor3.1}
The system (\ref{eq2}) has a unique minimal solution iff $$(\delta_{1}(\hat{x}), \delta_{2}(\hat{x}), \cdots, \delta_{n}(\hat{x})) \in S(A,b).$$ Moreover, $(\delta_{1}(\hat{x}), \delta_{2}(\hat{x}), \cdots, \delta_{n}(\hat{x}))$ is the unique minimal solution.
\end{corollary}
\begin{proof}
Let $y_j\in F_{j}(\hat{x})$. Then $(1,\cdots, 1, y_j, 1,\cdots,1)\in S(A,b)$. So that if $x =(x_1, x_2, \cdots, x_n)$ is the unique minimal solution of system \eqref{eq2}, then $x_j \le y_j$, thus by the arbitrariness of $y_j$, $x_j \le \delta_{j}(\hat{x})$. On the other hand, from Theorem \ref{th3.2}, $x_j \in F_{j}(\hat{x})$. Therefore, $x_j = \delta_{j}(\hat{x})$ for any $j \in J$.

Conversely, if $(\delta_{1}(\hat{x}), \delta_{2}(\hat{x}), \cdots, \delta_{n}(\hat{x})) \in S(A,b)$, then from Theorem \ref{th3.2}, it is straightforward that $(\delta_{1}(\hat{x}), \delta_{2}(\hat{x}), \cdots, \delta_{n}(\hat{x}))$ is the unique minimal solution.
\end{proof}

\section{Algorithm for calculating minimal solutions of system (\ref{eq2})}
This section first proves that for every $x\in S(A,b)$ there exists an $x_*\in \check{S}(A,b)$ satisfying $x_*\le x$, and then presents an algorithm to search for minimal solutions of a given one to the system (\ref{eq2}) with computational complexity $O(m^3)$ or $O(n^3)$.

 From Proposition \ref{pro3.1}, if $a_{it}>0$ for all $i\in I, t\in J$, and $a_{ij}+y_j-1>0$ for all $i\in I$, then $y = (y_1, y_2, \cdots, y_n)$ defined by (\ref{eq03}) is a minimal solution when $S(A,b) \ne \emptyset$. However, if there is an $i\in I$ such that $a_{ij}+y_j-1\leq 0$ then Formula (\ref{eq03}) may be invalid for constructing a minimal solution as shown by the following example.
\begin{example}\label{ex2}
\emph{Consider the following fuzzy relational inequalities:}
\begin{equation*}
\left\{\begin{array}{ll}
T_L(0.5, x_1) + T_L(0.7, x_{2}) + T_L(0.4, x_3)\geq 1,\\
T_L(0.3, x_1) +T_L(0.5, x_{2})+ T_L(0.9, x_3)\geq 1.3,\\
T_L(0.8, x_1) +T_L(0.6, x_{2}) + T_L(0.7, x_3)\geq 1.6.
\end{array}\right.
\end{equation*}\end{example}
Obviously, $(1, 1, 1)\in S(A,b)$. From Proposition \ref{pro3.1}, $$x_1=\max\limits_{i\in \{1, 2, 3\}}\{0, 1+b_i-\sum_{j\in \{1, 2, 3\}}a_{ij}\}=\max\{0.4, 0.6, 0.5\}=0.6.$$ Since $a_{21}+0.6-1<0$, Proposition \ref{pro3.1} is not suitable for determining whether $x=(0.6, 1, 1)$ is a minimal solution. In fact, it is easy to see that $x=(0.6, 1, 1)$ is a solution of the fuzzy relational inequalities, and from Theorem \ref{th3.1}, $x=(0.6, 1, 1)$ is not a minimal solution.

Therefore, we have to investigate another method for finding minimal solutions of system (\ref{eq2}). First, we have the following theorem.
\begin{theorem}\label{th4.1}
If $x = (x_1, x_2, \cdots, x_n)\in S(A,b)$. Then $x\in \check{S}(A,b)$ iff $x_j=\delta_j(x)$ for any $j\in J$.
\end{theorem}
\begin{proof}
Let $x\in\check{S}(A,b)$. Then $x_j\in F_{j}(x)$ for any $j\in J$. Thus, $x_j\geq \delta_j(x)$ and for any $x_{*j}\in F_j(x)$,
\begin{eqnarray*} &&T_L(a_{ij},\min\{x_j,x_{*j}\})+ \sum_{k\in J\setminus\{j\}}T_L(a_{ik}, x_{k})
 \\&=&\min\{T_L(a_{ij}, x_{j}),T_L(a_{ij}, x_{*j}) \}+ \sum_{k\in J\setminus\{j\}}T_L(a_{ik}, x_{k})\\&=&\min\{T_L(a_{ij}, x_{j}) +\sum_{k\in J\setminus\{j\}}T_L(a_{ik}, x_{k}), T_L(a_{ij}, x_{*j}) +\sum_{k\in J\setminus\{j\}}T_L(a_{ik}, x_{k})\}\\&\geq&\min\{b_i, b_i\}\\&=&b_i
 \end{eqnarray*}
for all $i\in I$. Therefore, $\min\{x_j,x_{*j}\}\in F_j(x)$ for all $j\in J$. Furthermore, define $x'=(x'_1, x'_2, \cdots, x'_n)$ with $$x'_k = \left\{ \begin{array}{l}
\min\{x_j,x_{*j}\}, k = j,\\
x_{k}, k\ne j.
\end{array} \right.$$
Then clearly $x'\le x$ and $x'\in S(A,b)$. Since $x\in\check{S}(A, b)$, $x'=x$, i.e., $x_j\le x_{*j}$. Hence by the arbitrariness of $x_{*j}$, we have $x_j\le \delta_j(x)$. Consequently, $x_{j} = \delta_j(x)$ for every $j\in J$.

Conversely, suppose that $x_{j} = \delta_j(x)$ for every $j\in J$. Let $y= (y_1, y_2, \cdots, y_n)\in S(A,b)$ be such that $y\le x$. Then $y_j\le x_j=\delta_j(x)$ for every $j\in J$. By Proposition \ref{pro3.2}, $y_{j} \in F_{j}(x)$, which means that $y_j\ge \delta_j(x)$. Thus $y_j= x_j$ for every $j\in J$, i.e., $y= x$. Therefore, $x\in \check{S}(A,b)$.
\end{proof}

From the proof of Theorem \ref{th4.1}, the following statement is true.
\begin{proposition}\label{pro4.1}
If $x = (x_1, x_2, \cdots, x_n)\in S(A,b)$, and define $y =(y_1, y_2, \cdots, y_n)$ with $$y_k = \left\{ \begin{array}{l}
\delta_j(x), k = j,\\
x_{k}, k\ne j.
\end{array} \right.$$
\end{proposition}
Then $y\in S(A,b)$ and $y\le x$.

Let $x = (x_1, x_2, \cdots, x_n)\in S(A,b)$, and for any $j\in J$, denote
\begin{equation*}
x^{(12\cdots j)}=\left\{\begin{array}{ll}
(\delta_1(x), x_2, \cdots, x_n), j=1,\\ \\
(\delta_1(x), \cdots, \delta_j(x^{(12\cdots (j-1))}), x_{j+1}, \cdots, x_n), j\ge 2.
\end{array}\right.
\end{equation*}
We then have the theorem as below.

\begin{theorem}\label{th4.2}
If $x = (x_1, x_2, \cdots, x_n)\in S(A,b)$, then there is an $x^*\in \check{S}(A,b)$ satisfying $x^*\le x$
\end{theorem}
\begin{proof} Let $x = (x_1, x_2, \cdots, x_n)\in S(A,b)$. If $x\in \check{S}(A,b)$, then obviously, there is an $x^*=x \in \check{S}(A,b)$ satisfying $x^*\leq x$.

If $x\notin \check{S}(A,b)$, then according to Proposition \ref{pro4.1}, $x^{(1)}=(\delta_1(x), x_2, \cdots, x_n)\in S(A,b).$ Again by Proposition \ref{pro4.1}, we have $x^{(12)}=(\delta_1(x), \delta_2(x^{(1)}), x_3, \cdots, x_n)\in S(A,b)$, and so on, we finally get $x^{(12\cdots n)}=(\delta_1(x), \delta_2(x^{(1)}), \cdots, \delta_n(x^{(12\cdots (n-1))}))\in S(A,b)$ and $x^{(12\cdots n)}\le x$.

Suppose that $y=(y_1, y_2, \cdots, y_n)\in S(A,b)$ satisfies $y\le x^{(12\cdots n)}$. Then for any $j\in J$, $y_j\le \delta_j(x^{(12\cdots (j-1))})$. Let $x^1=(y_1, x_2, \cdots, x_n)$. Then $y\le x^1$. Thus from Proposition \ref{pro2.1}, $x^1\in S(A,b)$, i.e., $y_1\in F_1(x)$, which means $\delta_1(x)\le y_1$. Hence, $\delta_1(x)= y_1$ since $y\le x^{(12\cdots n)}$, i.e., $y=(\delta_1(x), y_2, \cdots, y_n)$. Let $x^2= (\delta_1(x), y_2, x_3,\cdots, x_n)$. Then $y\le x^2$. Thus from Proposition \ref{pro2.1}, $x^2\in S(A,b)$, i.e., $y_2\in F_2(x^{(1)})$, which means $\delta_2(x^{(1)})\le y_2$. Hence, $\delta_2(x^{(1)})= y_2$ since $y\le x^{(12\cdots n)}$, i.e., $y=(\delta_1(x), \delta_2(x^{(1)}), \cdots, y_n)\in S(A,b)$. Repeating the process as above, we get $y_j=\delta_j(x^{(12\cdots (j-1))})$, $j = 3,\cdots ,n$. Therefore, $y=(\delta_1(x), \delta_2(x^{(1)}), \cdots, \delta_n(x^{(12\cdots (n-1))}))=x^{(12\cdots n)}$,
i.e., $x^{(12\cdots n)} \in \check{S}(A,b)$ and $x^{(12\cdots n)}\le x$.
\end{proof}

Applying Proposition \ref{pro2.1} and Theorem \ref{th4.2}, the following is true.
\begin{theorem}\label{th4.3}
\begin{eqnarray*}
S(A,b) =\bigcup_{\check{x}\in \check{S}(A, b)}\{x\in [0,1]^n \mid \check{x}\le x \le (1, 1, \cdots, 1)\}.
\end{eqnarray*}
\end{theorem}
\begin{theorem}
If $S(A,b)\neq \emptyset$ and $a_{ij}>0$ for all $i\in I$, $j\in j$, then the system (\ref{eq2}) has a unique solution iff there is an $i\in I$ such that $\sum_{j\in J}a_{ij}= b_i$.
\end{theorem}
\begin{proof}
 If $S(A,b)\neq \emptyset$, then $(1,1,\cdots, 1)\in S(A,b)$. If the system (\ref{eq2}) has a unique solution, then from Theorem \ref{th4.3}, $(1,1,\cdots, 1)$ is also a unique minimal solution of system (\ref{eq2}). By Theorem \ref{th3.1}, there is an $i\in I$ such that $\sum_{j\in J}T_L(a_{ij}, 1)= b_i$, i.e., $\sum_{j\in J}a_{ij}=b_i$.

Conversely, it is a straightforward matter from Theorem \ref{th2.2}.
\end{proof}

Let $x\in S(A,b)$. Then from the proof of Theorem \ref{th4.2}, we can summarize the following algorithm for calculating an $ \check{x}\in \check{S}(A, b)$ such that $ \check{x}\le x$.

\begin{algorithm}\label{alt1}Input $x = (x_1, x_2, \cdots, x_n)\in S(A,b)$. Output $ \check{x}$.\\
Step 1. Calculate $\tilde{x}= (\delta_{1}(\hat{x}), \delta_{2}(\hat{x}), \cdots, \delta_{n}(\hat{x}))$ defined by (\ref{eq4}). If $\tilde{x}\in S(A,b)$, then $ \check{x}:=\tilde{x}$, go to Step 8.\\
Step 2. j:=j+1\\
Step 3. Calculate \begin{eqnarray*}F_j(x^{(12\cdots (j-1))})&=& \{x_{*j} \in [0,1] |T_L(a_{i1}, \delta_1(x)) +\cdots+ T_L(a_{i(j-1)}, \delta_{j-1}(x^{(12\cdots (j-2))})\\&&+T_L(a_{ij}, x_{*j})+ \sum\limits_{k=j+1}^n T_L( a_{ik}, x_k) \ge b_i \mbox{ for any } i \in I\}\end{eqnarray*} where $x^{(0)}=x.$\\
Step 4. Calculate $$\delta_j(x^{(12\cdots (j-1))})=\min\{x_{*j}|x_{*j}\in F_j(x^{(12\cdots (j-1))})\}.$$\\
Step 5. $x^{(12\cdots j)}:=(\delta_1(x), \delta_2(x^{(1)}),\cdots, \delta_j(x^{(12\cdots (j-1))}, x_{j+1}, \cdots, x_n)$.\\
Step 6. Go to Step 2 when $j< n$.\\
Step 7. $ \check{x}:=x^{(12\cdots j)}$.\\
Step 8. Output $ \check{x}$.
\end{algorithm}
\begin{remark}\label{re1}
\emph{We can surely replace $(1,2,\cdots,n)$ in Algorithm \ref{alt1} by any of the permutations of $\{1,2,\cdots,n\}$, and just by repeating the steps as shown in Algorithm \ref{alt1}, we can obtain a minimal solution. However, using Algorithm \ref{alt1}, we can find at most $n!$ distinct minimal solutions since there are $n!$ different permutations and two different permutations may produce a same minimal solution.}
\end{remark}

The following two examples illustrate Algorithm \ref{alt1} and Remark  \ref{re1}, respectively.
\begin{example}
\emph{Consider the fuzzy relational inequalities in Example \ref{ex1}.}
\end{example}

It is obvious that $x=(0.8, 0.9, 1)$ is a solution of the fuzzy relational inequalities. With Algorithm \ref{alt1}, we can find a minimal solution as follows:\\
Step 1. Calculate $\delta_{1}(\hat{x})= 0.6$, $\delta_{2}(\hat{x})= 0.6$, $\delta_{3}(\hat{x})= 0.6$. Obviously, $(0.6, 0.6, 0.6)$ is not a solution.\\
Step 2. $j:=0+1$.\\
Step 3. Calculate $F_1(x)=\{x_1\mid 0.7\leq x_1\}. $\\
Step 4. Calculate $\delta_1(x)=0.7$.\\
Step 5. $x^{(1)}:=(0.7, 0.9, 1)$.\\
Step 6. $j=1< 3$, $j:=1+1$.\\
Step 7. Calculate $F_2(x^{(1)})=\{x_2\mid 0.9\leq x_2\}. $\\
Step 8. Calculate $\delta_2(x^{(1)})=0.9$.\\
Step 9. $x^{(12)}:=(0.7, 0.9, 1)$.\\
Step 10. $j=2< 3$, $j:=2+1$.\\
Step 11. Calculate $F_3(x^{(12)})=\{x_3\mid x_3=1\}. $\\
Step 12. Calculate $\delta_3(x^{(12)})=1$.\\
Step 13. $x^{(123)}:=(0.7, 0.9, 1)$.\\
Step 14. $j=3\nless 3$.\\
Step 15. $ \check{x}:=x^{(123)}$.\\
Step 16. Output $\check{x}=(0.7, 0.9, 1)$.\\
Thus, $\check{x}$ is a minimal solution of the fuzzy relational inequalities satisfying $\check{x}\leq x$.

\begin{example}\label{ex3}
\emph{Consider the fuzzy relational inequalities in Example \ref{ex2}.}
\end{example}

Obviously $x=(0.9, 0.9, 0.9)$ is a solution. With Algorithm \ref{alt1}, we can get three different minimal solutions.

Compute $\delta_{1}(\hat{x})= 0.5$, $\delta_{2}(\hat{x})= 0.6$, $\delta_{3}(\hat{x})= 0.6$. Obviously, $(0.5, 0.6, 0.6)$ is not a solution. Compute $$x^{(1)}=(\delta_{1}({x}), 0.9, 0.9)=(0.8, 0.9, 0.9),$$
$$x^{(12)}=(0.8, \delta_{2}({x^{(1)}}), 0.9)=(0.8, 0.9, 0.9),$$
and $$x^{(123)}=(0.8, 0.9, \delta_{3}({x^{(12)}}))=(0.8, 0.9, 0.9).$$

Similarly, we can also have the following minimal solutions:
 $$x^{(213)}=x^{(231)}=(0.9, 0.8, 0.9),$$
 $$x^{(321)}=x^{(312)}=(0.9, 0.9, 0.8),$$
 $$x^{(132)}=x^{(123)}=(0.8, 0.9, 0.9).$$
Therefore, $x^{(213)}$, $x^{(321)}$ and $x^{(132)}$ are three minimal solutions of the fuzzy relational inequalities fulfilling $x^{(213)}\leq x$, $x^{(321)}\leq x$ and $x^{(132)}\leq x$.

\begin{theorem}
 Algorithm \ref{alt1} terminates after $O(m^3)$ or $O(n^3)$ operations.
\end{theorem}
\begin{proof}
  The computational amount of Step 1 is $(4mn+m)n +4mn$. For every $j\in J$, the computational amount of Steps 3 and 4 is $4mn+m$. Then from Steps 2 to 6, the computational amount is $(4mn+m)n$. Thus the computational amount of Algorithm \ref{alt1} is $(4mn+m)n +4mn+(4mn+m)n= 8mn^2+6mn$. Thus, the computational complexity of Algorithm \ref{alt1} is $O(m^3)$ or $O(n^3)$.
\end{proof}
\section{Minimax programming problems with constraints of addition-{\L}ukasiewicz fuzzy relational inequalities}
In this section, we first transform the problem \eqref{eq6} into a single variable optimization problem, and show that every single variable optimization problem has a unique optimal solution. We then provide two algorithms to solve the single variable optimization problem, and obtain the greatest optimal solution of problem \eqref{eq6}. We further apply Algorithm \ref{alt1} for searching minimal optimal solutions of problem \eqref{eq6} which are less than or equal to the greatest optimal solution.

Denote the set of all the optimal solutions of problem \eqref{eq6} by $S^*(A,b)$. From Theorem \ref{th4.3}, one can see that $S(A,b)$ is a bounded set and the objective function of problem \eqref{eq6} is $Z(x)= x_1\vee x_2\vee \cdots\vee x_n$. Therefore, the optimal solution of problem \eqref{eq6} must exist, i.e., $S^*(A,b)\neq \emptyset$ if $S(A,b)\neq\emptyset$.

We first have the following three statements.
\begin{theorem}\label{the5.1}
If $S(A,b)\neq \emptyset$. Then there exists a minimal solution $x$ of system \eqref{eq2} such that $x\in S^*(A,b)$.
\end{theorem}
\begin{proof}
Take an arbitrary $x= (x_1, x_2, \cdots, x_n)\in S^*(A,b)$. If $x\in \check{S}(A, b)$, the proof
is completed. Otherwise, $x\notin \check{S}(A, b)$. Then by Theorem \ref{th4.2}, there exists $\check{x}\in \check{S}(A, b)$ such that  $\check{x}\leq x$. Hence, $$Z(\check{x})= \check{x}_1\vee \check{x}_2\vee \cdots\vee \check{x}_n\leq x_1\vee x_2\vee \cdots\vee x_n= Z(x).$$
If $Z(\check{x})= Z(x)$,  then $\check{x}\in S^*(A,b)$ and the proof is completed. Otherwise, $Z(\check{x})< Z(x)$. Because $\check{x}$ is a feasible solution and $x$ is an optimal solution of problem \eqref{eq6}, we get $Z(x)\leq Z(\check{x})$. This causes conflicts and the proof is complete.
\end{proof}

\begin{theorem}\label{the5.201}
Let $v\in S^*(A,b)$. If $x\in \check{S}(A,b)$ and $x \leq v$, then $x$ is a minimal optimal solution of problem \eqref{eq6}.
\end{theorem}
\begin{proof}
If $x\in \check{S}(A,b)$ and $x \leq v$, then $Z(x)\leq Z(v)$. Since $v\in S^*(A,b)$, $Z(v)\leq Z(x)$. Thus $Z(x)=Z(v)$.
Therefore, $x$ is a minimal optimal solution of problem \eqref{eq6} since $x\in \check{S}(A,b)$.
\end{proof}

From Theorem \ref{the5.1}, it is obvious that if the system \eqref{eq2} has a unique minimal solution $\check{x}$ then $\check{x}$ is the unique minimal optimal solution of problem \eqref{eq6}. Generally, one of the optimal solutions of problem \eqref{eq6} can be always selected from minimal solutions of system \eqref{eq2}.
 \begin{proposition}\label{pr5.1}
Supposing $S(A,b)\neq \emptyset$.
\begin{enumerate}
\item [(i)] Let $x=(x_1, x_2, \cdots, x_n)\in S^*(A,b)$. Then $x$ is the greatest optimal solution of problem \eqref{eq6} iff $x_1=x_2=\cdots=x_n=Z(x)$.
\item [(ii)] $S^*(A,b)$ has the greatest element, i.e., there exists $x\in S^*(A,b)$ such that $x$ is the greatest optimal solution of problem \eqref{eq6}.
 \end{enumerate}
\end{proposition}
\begin{proof}
(i) Let $x=(x_1, x_2, \cdots, x_n)\in S^*(A,b)$. Then $Z(x)=x_1\vee x_2\vee \cdots\vee x_n$, i.e., $x_j\leq Z(x)$ for any $j\in J$. If there is a $j\in J$ such that $x_j< Z(x)$, then let $x'=(x'_1,, x'_2, \cdots, x'_n)$ with $$x'_k= \left\{\begin{array}{ll}
Z(x), &k=j,\\
x_k, &k\in J\setminus \{j\}.
\end{array}\right.$$ Thus $x <x'$ and $x'\in S^*(A,b)$ since $Z(x')=Z(x)$. This contradicts the fact that $x$ is the greatest optimal solution. Therefore $x_1=x_2=\cdots=x_n=Z(x)$.

Conversely, if $x=(x_1, x_2, \cdots, x_n)\in S^*(A,b)$ with $x_j=Z(x)$ for any $j\in J$, then $Z(x)=y_1\vee y_2\vee \cdots\vee y_n$ for any $y=(y_1, y_2, \cdots, y_n)\in S^*(A,b)$. Therefore, $x_j\geq y_j$ for any $j\in J$, i.e., $x$ is the greatest optimal solution of problem \eqref{eq6}

(ii) From Theorem \ref{the5.1}, $S^*(A,b)\neq \emptyset$. For any $x=(x_1, x_2, \cdots, x_n)\in S^*(A,b)$, let $\bar{x}=(\bar{x}_1, \bar{x}_2, \cdots, \bar{x}_n)$ with $\bar{x}_j=Z(x)$ for any $j\in J$. Obviously, $\bar{x}\geq x$. Then by Proposition \ref{pro2.1}, $\bar{x}\in S(A,b)$, and this follows that $\bar{x}=(\bar{x}_1, \bar{x}_2, \cdots, \bar{x}_n)\in S^*(A,b)$ since $Z(\bar{x})= Z(x)$. Further, by (i) $\bar{x}$ is the greatest optimal solution of problem \eqref{eq6}.
\end{proof}

In what follows, we consider how to solve the problem \eqref{eq6} when $S(A,b)\neq \emptyset$.

The problem \eqref{eq6} can first be converted into a single variable optimization problem as below.
\begin{align}\label{eq7}
\mbox{Minimize}\quad&Z(\mathbf{y})= y \nonumber\\
\mbox{subject to}\quad& \left\{\begin{array}{ll}
T_L(a_{11}, y) + T_L(a_{12}, y) + \cdots + T_L(a_{1n}, y)\geq b_1,\\
T_L(a_{21}, y) +T_L(a_{22}, y) + \cdots + T_L(a_{2n}, y)\geq b_2,\\
\cdots\\
T_L(a_{m1}, y) +T_L(a_{m2}, y) + \cdots + T_L(a_{mn}, y)\geq b_m,
\end{array}\right.
\end{align}
where $\mathbf{y}=(y,y,\cdots, y)$ with $y\in [0,1]$. Furthermore, the problem \eqref{eq7} can be decomposed into $m$ minimax subproblems with the constraint of a single inequity as follows.
\begin{align}\label{eq8}
\mbox{Minimize}\quad&Z(\mathbf{y})= y \nonumber\\
\mbox{subject to}\quad& T_L(a_{i1}, y)+T_L(a_{i2}, y)+\cdots+T_L(a_{in}, y)\geq b_i,
\end{align}
where $\mathbf{y}=(y,y,\cdots, y)$ with $y\in [0,1]$ and $i\in I$. For any $i\in I$, denoted by $$S^i(A,b)=\{(y,y,\cdots, y)|T_L(a_{i1}, y)+T_L(a_{i2}, y)+\cdots+T_L(a_{in}, y)\geq b_i\}.$$

Now, we discuss the existence and uniqueness of the optimal solution of problem \eqref{eq8} and  develop an algorithm to find the optimal solution of problem \eqref{eq6}, with an illustrative numerical example.
\begin{proposition}\label{pr5.2}
 Problem \eqref{eq8} has a unique optimal solution iff $S^i(A,b)\neq \emptyset$.
\end{proposition}
\begin{proof}
 Suppose the problem \eqref{eq8} has a unique optimal solution $\mathbf{x}$. Then clearly $\mathbf{x}\in S^i(A,b)$.

Conversely, if $S^i(A,b)\neq \emptyset$, then by Theorem \ref{the5.1}, there is an $x\in S^i(A,b)$ such that $x$ is the optimal solution of problem \eqref{eq8}. Now, suppose both $\mathbf{v}=(v,v,\cdots, v)$ and $\mathbf{w}=(w,w,\cdots, w)$ are the optimal solutions of problem \eqref{eq8}. Then $z(\mathbf{v})= v\leq w= z(\mathbf{w})$ and $z(\mathbf{w}) = w\leq v= z(\mathbf{v})$. Hence, $\mathbf{w}=\mathbf{v}$. This verifies the uniqueness of the optimal solution.
\end{proof}

In the following, we always suppose that $S^i(A,b)\neq \emptyset$ with $i\in I$.
Let $\mathbf{u}= (u, u, \cdots, u)\in [0, 1]^n$ with
\begin{equation}\label{eq9}
 u = \frac{b_i-\sum_{j\in J} a_{ij}}{n}+1,
\end{equation}
and denote \begin{equation}\label{eq10}J(\mathbf{u})=\{j\in J|u\geq 1-a_{ij}\}.\end{equation}  By the definition of $J(\mathbf{u})$, we know that  $T_L(a_{ij}, u)=\max\{a_{ij}+ u-1,0\}=0$ for any $j\in J\setminus J(\mathbf{u})$ and $J(\mathbf{y})\subseteq J(\mathbf{u})$ for any $\mathbf{y}< \mathbf{u}$. Moreover, since $b_i>0$, $a_{ij}\in [0,1]$ and $\sum _{j\in J} a_{ij}\geq b_i$, $-1\leq \frac{b_i-\sum _{j\in J} a_{ij}}{n}\leq 0$. Then $u\in [0, 1]$.

\begin{theorem}\label{th5.1}
The following three statements hold:
\begin{enumerate}\item [(i)]$\mathbf{u}= (u, u, \cdots, u)\in S^i(A,b)$ and $J(\mathbf{u})\neq \emptyset$.
\item [(ii)] If $\sum_{j\in J(\mathbf{u})} T_L(a_{ij}, u)>b_i$, then there is a $\mathbf{u^1}=(u^1, u^1, \cdots, u^1)<\mathbf{u}$ such that $\mathbf{u^1}\in S^i(A,b)$.
\item [(iii)] $\mathbf{u}$ is the unique optimal solution of problem \eqref{eq8} iff $\sum_{j\in J(\mathbf{u})} T_L(a_{ij}, u)=b_i$.
\end{enumerate}\end{theorem}
\begin{proof}
\begin{enumerate}\item [(i)] From Formula \eqref{eq9},
\begin{eqnarray*}
&&T_L(a_{i1}, u)+T_L(a_{i2}, u)+\cdots+T_L(a_{in}, u)\\&=&\max\{a_{i1}+u-1,0\}+\max\{a_{i2}+u-1,0\}+\cdots+\max\{a_{in}+u-1,0\}\\&\geq&
a_{i1}+u-1+a_{i2}+u-1+\cdots+a_{in}+u-1\\&=&\sum\limits_{j\in J} a_{ij}+nu-n\\&=&\sum\limits_{j\in J} a_{ij}+n(\frac{b_i-\sum _{j\in J} a_{ij}}{n}+1)-n\\&=&b_i.
\end{eqnarray*}
Therefore,  $\mathbf{u}\in S^i(A,b)$. This follows that $J(\mathbf{u})\neq \emptyset$ since $b_i>0$.
\item [(ii)] Since $T_L(a_{ij}, u)=\max\{a_{ij}+u-1,0\}=a_{ij}+u-1$ for any $j\in J(\mathbf{u})$, we have $$\sum\limits_{j\in J(\mathbf{u})} T_L(a_{ij}, u)= \sum \limits_{j\in J(\mathbf{u})}(a_{ij}+u-1)=\sum\limits_{j\in J(\mathbf{u})}a_{ij}+|J(\mathbf{u})|u-|J(\mathbf{u})|>b_i.$$ Therefore, $u>\frac{b_i-\sum _{j\in J(\mathbf{u})} a_{ij}}{|J(\mathbf{u})|}+1$.

Let $\mathbf{u^1}=(u^1, u^1, \cdots, u^1)\in [0, 1]^n$ with $u^1=\frac{b_i-\sum_{j\in J(\mathbf{u})} a_{ij}}{|J(\mathbf{u})|}+1$.
Obviously, $\mathbf{u^1}<\mathbf{u}$. Because
\begin{eqnarray*}
\sum\limits_{j\in J} T_L(a_{ij}, u^1)&\geq&\sum\limits_{j\in J(\mathbf{u})} T_L(a_{ij}, u^1)\\&=&\sum\limits_{j\in J(\mathbf{u})}\max\{a_{ij}+u^1-1,0\}\\&\geq& \sum\limits_{j\in J(\mathbf{u})}(a_{i1}+u^1-1)\\&=&\sum\limits_{j\in J(\mathbf{u})} a_{ij}+|J(\mathbf{u})|u^1-|J(\mathbf{u})|\\&=&\sum\limits_{j\in J(\mathbf{u})} a_{ij}+|J(\mathbf{u})|(\frac{b_i-\sum_{j\in J(\mathbf{u})} a_{ij}}{|J(\mathbf{u})|}+1)-|J(\mathbf{u})|\\&=&b_i,
\end{eqnarray*}
we have $\mathbf{u^1}\in S^i(A, b)$.

\item [(iii)] By (i), we know $\mathbf{u}\in S^i(A,b)$. If $\mathbf{u}$ is the unique optimal solution of problem \eqref{eq8} and $\sum_{j\in J(\mathbf{u})} T_L(a_{ij}, u)>b_i$, then by (ii), there is a $\mathbf{u^1}=(u^1, u^1, \cdots, u^1)<\mathbf{u}$ such that $\mathbf{u^1}\in S^i(A,b)$. This follows that $Z(\mathbf{u^1})< Z(\mathbf{u})$, contrary to the fact that $\mathbf{u}$ is the unique optimal solution of problem \eqref{eq8}. Therefore, $\sum_{j\in J(\mathbf{u})} T_L(a_{ij}, u)=b_i$

Conversely, it is clearly that $\mathbf{u}\in S^i(A,b)$ since $\sum_{j\in J(\mathbf{u})} T_L(a_{ij}, u)=b_i$. Supposing that $\mathbf{y}=(y, y, \cdots, y)\in S^i(A,b)$ with $y<u$. By Formula \eqref{eq10}, we can get $J(\mathbf{y})\subseteq J(\mathbf{u})$. Since $\sum _{j\in J(\mathbf{u})} T_L(a_{ij}, u)=b_i$, $$\sum \limits_{j\in J}T_L(a_{ij}, y)=\sum \limits_{j\in J(\mathbf{y})}T_L(a_{ij}, y)<\sum\limits_{j\in J(\mathbf{y})}T_L(a_{ij}, u)\leq\sum \limits_{j\in J(\mathbf{u})}T_L(a_{ij}, u)=b_i,$$ contrary to $\mathbf{y}\in S^i(A,b)$.
Thus $\mathbf{u}$ is an optimal solution of problem \eqref{eq8}. Therefore, by Proposition \ref{pr5.2} $\mathbf{u}$ is the unique optimal solution of problem \eqref{eq8}.
\end{enumerate}\end{proof}
\begin{theorem}\label{th5.2}
Let $\mathbf{u^1}=(u^1, u^1, \cdots, u^1)\in [0, 1]^n$ with $$u^1=\frac{b_i-\sum_{j\in J(\mathbf{u})} a_{ij}}{|J(\mathbf{u})|}+1.$$ Then $\mathbf{u^1}$ is the unique optimal solution of problem \eqref{eq8} iff $J(\mathbf{u^1})=J(\mathbf{u})$.
\end{theorem}
\begin{proof}
By Theorem \ref{th5.1} (i), we know that $\mathbf{u}\in S^i(A,b)$. Let $\mathbf{u^1}=(u^1, u^1, \cdots, u^1)$ with $u^1=\frac{b_i-\sum_{j\in J(\mathbf{u})} a_{ij}}{|J(\mathbf{u})|}+1$. According to the proof of Theorem \ref{th5.1} (ii), we have $\mathbf{u^1}\leq \mathbf{u}$ and $\mathbf{u^1}\in S^i(A,b)$. Further, by Formula \eqref{eq10}, $J(\mathbf{u^1})\subseteq J(\mathbf{u})$.

If $J(\mathbf{u^1})=J(\mathbf{u})$, then
\begin{eqnarray*}\label{eq11}
\sum_{j\in J(\mathbf{u^1})} T_L(a_{ij}, u^1)&=&\sum_{j\in J(\mathbf{u})}\max\{a_{ij}+u^1-1,0\}\\&=&\sum_{j\in J(\mathbf{u})}(a_{ij}+u^1-1)\\&=&\sum_{j\in J(\mathbf{u})}a_{ij}+|J(\mathbf{u})|u^1-|J(\mathbf{u})|\\&=&\sum_{j\in J(\mathbf{u})}a_{ij}+|J(\mathbf{u})|(\frac{b_i-\sum_{j\in J(\mathbf{u})} a_{ij}}{|J(\mathbf{u})|}+1)-|J(\mathbf{u})|\\&=&b_i.
 \end{eqnarray*}
Then by Theorem \ref{th5.1} (iii), $\mathbf{u^1}$ is the unique optimal solution of problem \eqref{eq8}.

Conversely, suppose that $\mathbf{u^1}$ is the unique optimal solution of problem \eqref{eq8}. Then by Theorem \ref{th5.1} (iii), $\sum_{j\in J(\mathbf{u^1})} T_L(a_{ij}, u^1)=b_i$. Assume that there exists a $j\in J(\mathbf{u})\setminus J(\mathbf{u^1})$. Then $a_{ij}+u^1-1< 0$, i.e., $$\sum_{j\in J(\mathbf{u^1})} T_L(a_{ij}, u^1)=\sum_{j\in J(\mathbf{u^1})} (a_{ij}+u^1-1)>\sum_{j\in J(\mathbf{u})}(a_{ij}+u^1-1)=b_i,$$ a contradiction. Therefore, $J(\mathbf{u})\setminus J(\mathbf{u^1})=\emptyset$, i.e., $J(\mathbf{u^1})=J(\mathbf{u})$.
\end{proof}

\begin{theorem}\label{th5.3}
 There is a $\mathbf{y}$ with $\mathbf{y}\leq \mathbf{u}$ such that $\mathbf{y}$ is the unique optimal solution of problem \eqref{eq8}.
\end{theorem}
\begin{proof}
  From Theorem \ref{th5.1} (i), $\mathbf{u}\in S^i(A,b)$. Then $\sum_{j\in J(\mathbf{u})} T_L(a_{ij}, u)\geq b_i$. If $\sum_{j\in J(\mathbf{u})} T_L(a_{ij}, u)=b_i$, then by Theorem \ref{th5.1} (iii) $\mathbf{y}= \mathbf{u}$ is the unique optimal solution of problem \eqref{eq8}. Otherwise, $\sum_{j\in J(\mathbf{u})} T_L(a_{ij}, u)>b_i$. By Theorem \ref{th5.1} (ii) there is a $\mathbf{u^1}=(u^1, u^1, \cdots, u^1)$ such that $\mathbf{u^1}<\mathbf{u}$ and $\mathbf{u^1}\in S^i(A,b)$. From Theorem \ref{th5.2}, if $J(\mathbf{u^1})=J(\mathbf{u})$ then $\mathbf{u^1}$ is the unique optimal solution of problem \eqref{eq8}. Otherwise, $\sum_{j\in J(\mathbf{u^1})} T_L(a_{ij}, u^1)>b_i$ and $J(\mathbf{u^1})\subset J(\mathbf{u})$. Thus by Theorem \ref{th5.1} (ii), there exists a $\mathbf{u}^2=(u^2, u^2, \cdots, u^2)$ such that $\mathbf{u^2}<\mathbf{u^1}$ and $\mathbf{u^2}\in S^i(A,b)$. Again from Theorem \ref{th5.2}, if $J(\mathbf{u^2})=J(\mathbf{u^1})$ then $\mathbf{u^2}$ is the unique optimal solution of problem \eqref{eq8}. Otherwise, $\sum_{j\in J(\mathbf{u^2})} T_L(a_{ij}, u^2)>b_i$ and $J(\mathbf{u^2})\subset J(\mathbf{u^1})$. Continuing the process as above, there must exist a $\mathbf{u^k}$  such that $$J(\mathbf{u^k})=J(\mathbf{u^{k-1}})\subset J(\mathbf{u^{k-2}})\subset\cdots\subset J(\mathbf{u^1})\subset J(\mathbf{u})$$ since $\cdots\subset J(\mathbf{u^2})\subset J(\mathbf{u^1})\subset J(\mathbf{u}) \subseteq J$ and $J$ is a finite set. Therefore, by Theorem \ref{th5.2} $\mathbf{u^k}$ is the unique optimal solution of problem \eqref{eq8}. Let $\mathbf{y}=\mathbf{u^k}$. Then $\mathbf{y}\leq \mathbf{u}$ and $\mathbf{y}$ is the unique optimal solution of problem \eqref{eq8}.
\end{proof}

Based on the proof of Theorem \ref{th5.3}, we can summarize the following algorithm for searching the unique optimal solution of problem \eqref{eq8}.
\begin{algorithm}\label{al2}
Input $a_{i1}, a_{i2}, \cdots, a_{in}$ and $b_i$ with $i\in I$. Output $\mathbf{u^k}$.\\
Step 1. Check the solvability of problem \eqref{eq8}. If $S^i(A,b)=\emptyset$, then stop.\\
Step 2. Compute $\mathbf{u^0}=(u^0, u^0, \cdots, u^0)$ with $u^0=u$ defined by (\ref{eq9}), $k: = 0.$ \\
Step 3. Compute $J(\mathbf{u^k})$ defined by (\ref{eq10}). If $\sum_{j\in J(\mathbf{u^k})}T_L(a_{ij}, u^k)=b_i$, then go to Step 6.\\
Step 4. $k := k + 1$, compute $\mathbf{u^{k}}=(u^k, u^k, \cdots, u^k)$ with $u^{k}=\frac{b_i-\sum_{j\in J(\mathbf{u^{k-1}})} a_{ij}}{|J(\mathbf{u^{k-1}})|}+1$. \\
Step 5. Compute $J(\mathbf{u^{k}})$ defined by (\ref{eq10}). If $J(\mathbf{u^{k}})\subset J(\mathbf{u^{k-1}})$, then go to Step 4.\\
Step 6. Output $\mathbf{u^k}$.
\end{algorithm}
\begin{theorem}\label{thm5.5}
 Algorithm \ref{al2} terminates after $O(n^2)$ operations.
\end{theorem}
\begin{proof}
Checking $S^i(A,b)\neq \emptyset$ in Step 1 costs $4n$ operations. In Step 2, computing $\mathbf{u^0}$ needs $n+2$ operations. In Step 3, computing $J(\mathbf{u^k})$ and checking whether $\sum_{j\in J(\mathbf{u^k})}T_L(a_{ij}, u^k)=b_i$ hold or not, we take $6n$ operations. The computational amount of Steps 4 and 5 is $(3n+3)\times n$ since the number of loops in Steps 4 and 5 is less than or equal to $n$ by Theorem \ref{th5.3}. Thus, obtaining the optimal solution $\mathbf{u^k}$ requires
$$4n+n+2+6n+(3n+3)n=3n^2+14n+2$$ operations.
Therefore, the computational complexity of Algorithm \ref{al2} is $O(n^2)$.
\end{proof}

\begin{example}\label{ex4}
\emph{Consider the single variable optimization problem as follows:}
\begin{align*}
\emph{Minimize}\quad&Z(\mathbf{y})= y\\
\emph{subject to}\quad&\begin{array}{l}
T_L(0, y) + T_L(0.1, y) + T_L(0.5, y)+ T_L(0.8, y)+ T_L(0.6, y)+ T_L(0.3, y)\geq 1.4,
\end{array}
\end{align*}
where $y\in [0, 1]$.
\end{example}
Step 1. Obviously, $(1, 1, 1, 1, 1, 1)\in S^i(A,b)$.\\
Step 2. $u^0=\frac{b_i-\sum_{j\in J} a_{ij}}{n}+1=\frac{1.4-0-0.1-0.5-0.8-0.6-0.3}{6}+1=0.85$, i.e., $$\mathbf{u^0}=(0.85,0.85,0.85,0.85,0.85, 0.85).$$\\
Step 3. $J(\mathbf{u^0})=\{j\in J|u^0\geq 1-a_{ij}\}=\{3, 4, 5, 6\}$, and $\sum_{j\in J(\mathbf{u^0})}T_L(a_{ij}, u^0)=1.6>1.4$.\\
Step 4. $k=1$, $u^{1}=\frac{b_i-\sum_{j\in J(\mathbf{u^{0}})} a_{ij}}{|J(\mathbf{u^{0}})|}+1=\frac{1.4-0.5-0.8-0.6-0.3}{4}+1=0.8$, i.e., $$\mathbf{u^{1}}=(0.8, 0.8, 0.8, 0.8, 0.8, 0.8).$$
Step 5. $J(\mathbf{u^1})=\{j\in J|u^0\geq 1-a_{ij}\}=\{3, 4, 5, 6\}$, and $J(\mathbf{u^1})=J(\mathbf{u^0})$.\\
Step 6. Output $\mathbf{u^1}=(0.8, 0.8, 0.8, 0.8, 0.8, 0.8)$.

For solving the single variable problem \eqref{eq7}, we need the following proposition.
\begin{proposition}\label{pr5.3}
 Problem \eqref{eq7} has a unique optimal solution iff $S(A,b)\neq \emptyset$.
\end{proposition}
\begin{proof}
 This is shown in complete analogy to the proof of Proposition \ref{pr5.2}.
\end{proof}

If $S(A,b)\neq \emptyset$, then $S^i(A,b)\neq \emptyset$ for every $i\in I$. Thus from Proposition \ref{pr5.2}, we can let $\mathbf{u^{i*}}=(u^{i*}, u^{i*}, \cdots, u^{i*})$ with $i\in I$ be the unique optimal solution of problem \eqref{eq8}. Denoted by $\mathbf{u^{*}}=(u^{*}, u^{*}, \cdots, u^{*})$ with $u^*=\bigvee_{i\in I}u^{i*}$. Then we have the following theorem.
\begin{theorem}\label{th5.5}
 Let $S(A,b)\neq \emptyset$. Then $\mathbf{u^{*}}$ is the unique optimal solution of problem \eqref{eq7}. Further, $\mathbf{u^{*}}$ is the greatest optimal solution of problem \eqref{eq6}.
\end{theorem}
\begin{proof}
Obviously, $\mathbf{u^{*}}\geq \mathbf{u^{i*}}$ for any $i\in I$. Since $\mathbf{u^{i*}}$ is the unique optimal solution of problem \eqref{eq8}, we have $$T_L(a_{i1}, u^*) +T_L(a_{i2}, u^*) + \cdots + T_L(a_{in}, u^*)\geq T_L(a_{i1}, u^{i*}) +T_L(a_{i2}, u^{i*}) + \cdots + T_L(a_{in}, u^{i*})\geq b_i$$ for any $i\in I$. Hence, $\mathbf{u^{*}}$ is a feasible solution of problem \eqref{eq7}.

Next, we show that $\mathbf{u^{*}}$ is the unique optimal solution of problem \eqref{eq7}.

 Suppose that $\mathbf{y}=(y, y, \cdots, y)$ is an arbitrary feasible solution of problem \eqref{eq7}. Then $$T_L(a_{i1}, y) +T_L(a_{i2}, y) + \cdots + T_L(a_{in}, y)\geq b_i \mbox{~for any~} i\in I.$$
 Thus, $\mathbf{y}$ is a feasible solution of problem \eqref{eq8} for all $i\in I$. This follows that $Z(\mathbf{y}) \geq Z(\mathbf{u^{i*}})$ for any $i\in I$ since $\mathbf{u^{i*}}$ is the unique optimal solution of problem \eqref{eq8}. Consequently, $y\geq \bigvee_{i\in I}u^{i*}=u^*$, and this means that $Z(\mathbf{y}) \geq Z(\mathbf{u^{*}})$. Therefore, by Proposition \ref{pr5.3}, $\mathbf{u^{*}}$ is the unique optimal solution of problem \eqref{eq7}.

 Now, let $x=(x_1, x_2,\cdots, x_n)\in S^*(A,b)$. Then by Proposition \ref{pr5.1}, the greatest optimal solution of problem \eqref{eq6} is $\mathbf{y}=(y, y, \cdots, y)$ with $y=Z(x)$. Thus $Z(\mathbf{y})=Z(x)$, and $Z(\mathbf{y})\leq Z(\mathbf{u^{*}})$ since $\mathbf{u^*}$ is a feasible solution of problem \eqref{eq6}. The inequality implies that $y\leq u^*$, which follows that $\mathbf{y}\leq \mathbf{u^{*}}$.
On the other hand, it is obvious that $\mathbf{y}$ is also a feasible solution of problem \eqref{eq7}. This follows that $\mathbf{u^{*}}\leq \mathbf{y}$ since $\mathbf{u^{*}}$ is the unique optimal solution of problem \eqref{eq7}. Consequently, $\mathbf{y}=\mathbf{u^{*}}$ and $\mathbf{u^{*}}$ is the greatest optimal solution of problem \eqref{eq6}.
\end{proof}

Based on the above results and Algorithm \ref{al2}, we establish the following algorithm for finding the unique optimal solution of problem \eqref{eq7} which is also the greatest optimal solution of problem \eqref{eq6}.
\begin{algorithm}\label{al3}
Input $A=(a_{ij})_{m \times n}$ and $b=(b_1, b_2, \cdots, b_m)$. Output $\mathbf{u^*}$.\\
Step 1: Check the  solvability of problem \eqref{eq6}. If $S(A,b)= \emptyset$, then stop.\\
Step 2: Convert the problem \eqref{eq6} into the problem \eqref{eq7}.\\
Step 3: Decompose the problem \eqref{eq7} into $m$ subproblems, i.e., the problem \eqref{eq8}.\\
Step 4: For every $i\in I$, compute the unique optimal solution $\mathbf{u^{i*}}=(u^{i*}, u^{i*}, \cdots, u^{i*})$ of problem \eqref{eq8} by Algorithm \ref{al2}.\\
Step 5: Compute $\mathbf{u^{*}}=(u^{*}, u^{*}, \cdots, u^{*})$ with $u^*=\bigvee_{i\in I}u^{i*}$.\\
Step 6: Output $\mathbf{u^*}$.
\end{algorithm}
\begin{theorem}
 Algorithm \ref{al3} terminates after $O(n^3)$ or $O(m^3)$ operations.
\end{theorem}
\begin{proof}
From Theorem \ref{thm5.5}, one can clarify that the total computational complexity of Algorithm \ref{al3} is $O(n^3)$ or $O(m^3)$ operations since Step 4 is the key process of Algorithm \ref{al3} and its computational amount is $mn^2$.
\end{proof}

The following two numerical examples illustrate Algorithm \ref{al3}.
\begin{example}\label{ex5}
\emph{Consider a P2P network system with
3 terminals described by the addition-{\L}ukasiewicz fuzzy relational inequalities as follows.
$$A \odot {x^T} \geq {b^T}$$
where $$A=\left (
\begin{array}{cccc}
  0.5 & 0.7 & 0.4\\
  0.3 & 0.5 & 0.9\\
  0.8& 0.6 & 0.7\\
\end{array}
\right),$$
$$b=(1, 1.3, 1.6).$$
Our goal is to find an optimal solution of the following minimax programming problem:}
\begin{align}\label{eq11}
\emph{Minimize}\quad & Z(\mathbf{x})= x_1\vee x_2\vee x_3\nonumber\\
\emph{subject to}\quad & A \odot {x^T} \geq {b^T}.
\end{align}
\emph{where $x=(x_1, x_2, x_3)\in [0, 1]^3$.}
\end{example}
Step 1. Obviously, $(1, 1, 1)\in S(A,b)$.\\
Step 2. Convert the problem \eqref{eq11} into the following version:
\begin{align}\label{eq12}
\mbox{Minimize}\quad&Z(\mathbf{y})=y\nonumber\\
\mbox{subject to}\quad&\left\{\begin{array}{ll}
T_L(0.5, y) + T_L(0.7, y) + T_L(0.4, y)\geq 1,\\
T_L(0.3, y) +T_L(0.5, y)+ T_L(0.9, y)\geq 1.3,\\
T_L(0.8, y) +T_L(0.6,y) + T_L(0.7, y)\geq 1.6,
\end{array}\right.
\end{align}where $y\in [0, 1]$.\\
Step 3. Decompose the problem \eqref{eq12} into 3 subproblems.\\
\begin{align*}
(1)~\mbox{Minimize}\quad&Z(\mathbf{y})=y\nonumber\\
\mbox{subject to}\quad& \left\{\begin{array}{ll}
T_L(0.5, y) + T_L(0.7, y) + T_L(0.4, y)\geq 1,\\
y\in [0, 1].
\end{array}\right.
\end{align*}
\begin{align*}
(2)~\mbox{Minimize}\quad&Z(\mathbf{y})=y\nonumber\\
\mbox{subject to}\quad& \left\{\begin{array}{ll}
T_L(0.3, y) +T_L(0.5, y)+ T_L(0.9, y)\geq 1.3,\\
y\in [0, 1].
\end{array}\right.
\end{align*}
\begin{align*}
(3)~\mbox{Minimize}\quad&Z(\mathbf{y})=y\nonumber\\
\mbox{subject to}\quad& \left\{\begin{array}{ll}
T_L(0.8, y) +T_L(0.6,y) + T_L(0.7, y)\geq 1.6,\\
y\in [0, 1].
\end{array}\right.
\end{align*}
Step 4. By Algorithm \ref{al2}, we can get the unique optimal solution $\mathbf{u^{i*}}=(u^{i*}, u^{i*}, \cdots, u^{i*})$ with $i\in\{1,2,3\}$:
$\mathbf{u^{1*}}=(0.8, 0.8, 0.8)$, $\mathbf{u^{2*}}=(\frac{13}{15}, \frac{13}{15}, \frac{13}{15})$, and $\mathbf{u^{3*}}=(\frac{5}{6}, \frac{5}{6}, \frac{5}{6})$.\\
Step 5. $u^*=u^{1*}\vee u^{2*}\vee u^{3*}=\frac{13}{15}$, i.e., $\mathbf{u^{*}}=(\frac{13}{15}, \frac{13}{15}, \frac{13}{15})$.\\
Step 6: Output $\mathbf{u^*}=(\frac{13}{15}, \frac{13}{15}, \frac{13}{15})$.

Therefore, $(\frac{13}{15}, \frac{13}{15}, \frac{13}{15})$ is the greatest optimal solution of \eqref{eq11}.
\begin{example}\label{ex6}
\emph{Consider the P2P network system with 5 users and its corresponding minimax programming problem with addition-{\L}ukasiewicz fuzzy relational inequality constraints as below.
\begin{align}\label{eq13}
\mbox{Minimize}\quad&Z(\mathbf{x})= x_1\vee x_2\vee x_3\vee x_4\vee x_5\nonumber\\
\mbox{subject to}\quad&A \odot {x^T} \geq {b^T},
\end{align}where $x=(x_1, x_2, x_3, x_4, x_5)\in [0, 1]^5$,
$$A=\left( {\begin{array}{*{20}{c}}
0&0.6&0.5&0,1&0.7\\
0.4&0.8&0.7&0.6&0.5\\
0.3&0.2&0.3&0.8&0.2\\
0.7&0.5&0.5&0.4&0.7
\end{array}} \right)$$ and
$$b=(1.3, 1.5, 0.8, 1.6).$$}
\end{example}
Step 1. Obviously, $(1, 1, 1, 1, 1)\in S(A,b)$.\\
Step 2. Convert the problem \eqref{eq13} into the following version:
\begin{align}\label{eq14}
\mbox{Minimize}\, &Z(\mathbf{y})=y\nonumber\\
\mbox{subject to}\, &\left\{\begin{array}{ll}
T_L(0, y) + T_L(0.6, y) + T_L(0.5, y)+ T_L(0.1, y)+ T_L(0.7, y)\geq 1.3,\\
T_L(0.4, y) +T_L(0.8, y)+ T_L(0.7, y)+ T_L(0.6, y)+ T_L(0.5, y)\geq 1.5,\\
T_L(0.3, y) +T_L(0.2, y) + T_L(0.3, y)+ T_L(0.8, y)+ T_L(0.2, y)\geq 0.8,\\
T_L(0.7, y) +T_L(0.5, y)+ T_L(0.5, y)+ T_L(0.4, y)+ T_L(0.7, y)\geq 1.6,
\end{array}\right.
\end{align}where $y\in [0, 1]$.\\
Step 3. Decompose the problem \eqref{eq14} into 4 subproblems.\\
\begin{align*} (1)~\mbox{Minimize}\quad&Z(\mathbf{y})=y\nonumber\\
\mbox{subject to}\quad& \left\{\begin{array}{ll}
T_L(0, y) + T_L(0.6, y) + T_L(0.5, y)+ T_L(0.1, y)+ T_L(0.7, y)\geq 1.3,\\
y\in [0, 1].
\end{array}\right.
\end{align*}
\begin{align*} (2)~\mbox{Minimize}\quad&Z(\mathbf{y})=y\nonumber\\
\mbox{subject to}\quad& \left\{\begin{array}{ll}
T_L(0.4, y) +T_L(0.8, y)+ T_L(0.7, y)+ T_L(0.6, y)+ T_L(0.5, y)\geq 1.5,\\
y\in [0, 1].
\end{array}\right.
\end{align*}
\begin{align*}(3)~
\mbox{Minimize}\quad&Z(\mathbf{y})=y\nonumber\\
\mbox{subject to}\quad& \left\{\begin{array}{ll}
T_L(0.3, y) +T_L(0.2, y) + T_L(0.3, y)+ T_L(0.8, y)+ T_L(0.2, y)\geq 0.8,\\
y\in [0, 1].
\end{array}\right.
\end{align*}
\begin{align*}(4)~
\mbox{Minimize}\quad&Z(\mathbf{y})=y\nonumber\\
\mbox{subject to}\quad& \left\{\begin{array}{ll}
T_L(0.7, y) +T_L(0.5, y)+ T_L(0.5, y)+ T_L(0.4, y)+ T_L(0.7, y)\geq 1.6,\\
y\in [0, 1].
\end{array}\right.
\end{align*}
Step 4. By Algorithm \ref{al2}, we can get the unique optimal solution $\mathbf{u^{i*}}=(u^{i*}, u^{i*}, \cdots, u^{i*})$ with $i\in\{1,2,3,4\}$:
$\mathbf{u^{1*}}=(\frac{5}{6}, \frac{5}{6}, \frac{5}{6}, \frac{5}{6}, \frac{5}{6})$, $\mathbf{u^{2*}}=(0.7, 0.7, 0.7, 0.7, 0.7)$, \\ $\mathbf{u^{3*}}=(0.8, 0.8, 0.8, 0.8, 0.8)$, and $\mathbf{u^{4*}}=(0.76, 0.76, 0.76, 0.76, 0.76)$.\\
Step 5. $u^*=u^{1*}\vee u^{2*}\vee u^{3*}\vee u^{4*}=\frac{5}{6}$, i.e., $\mathbf{u^{*}}=(\frac{5}{6}, \frac{5}{6}, \frac{5}{6}, \frac{5}{6}, \frac{5}{6})$.\\
Step 6: Output $\mathbf{u^*}=(\frac{5}{6}, \frac{5}{6}, \frac{5}{6}, \frac{5}{6}, \frac{5}{6})$.

Therefore, $(\frac{5}{6}, \frac{5}{6}, \frac{5}{6}, \frac{5}{6}, \frac{5}{6})$ is the greatest optimal solution of \eqref{eq13}.

From Theorem \ref{th4.2}, there exists an $x\in \check{S}(A,b)$ such that $x\le \mathbf{u^*}$ since $\mathbf{u^*}\in S^*(A,b)$. Moreover, by Theorem \ref{the5.201} $x$ is a minimal optimal solution of problem \eqref{eq6}. In particular, from Theorems \ref{the5.1}, \ref{the5.201} and \ref{th5.5} we have the following theorem.
\begin{theorem}\label{th5.7}
$\mathbf{u^{*}}$ is the unique optimal solution of problem \eqref{eq6} iff $\mathbf{u^{*}}\in \check{S}(A,b)$.
\end{theorem}

\begin{remark}
\begin{enumerate}
\item [(i)]\emph{In Example \ref{ex5}, from Theorem \ref{th5.7} $\mathbf{u^{*}}$ is the unique optimal solution of problem \eqref{eq11} since $\mathbf{u^*}=(\frac{13}{15}, \frac{13}{15}, \frac{13}{15})$ is a minimal solution of system \eqref{eq2} by Theorem \ref{th3.1}. In Example \ref{ex6}, $(\frac{5}{6}, \frac{5}{6}, \frac{5}{6}, \frac{5}{6}, \frac{5}{6})$ is an optimal solution of problem \eqref{eq13}, but it is not a minimal optimal solution of problem \eqref{eq13}.}
\item [(ii)] \emph{By Algorithm \ref{al3}, we can first compute the greatest optimal solution $\mathbf{u^*}$ of problem \eqref{eq6}. Then applying Algorithm \ref{alt1}, we can obtain some $\check{x}\in \check{S}(A,b)$ with $\check{x}\le \mathbf{u^*}$, which are minimal optimal solutions of problem \eqref{eq6} by Theorem \ref{the5.201}.}
\end{enumerate}
\end{remark}
\begin{example}\label{ex7}
\emph{Consider some minimal optimal solutions of Example \ref{ex6}.}
\end{example}

By Algorithm \ref{al3}, we obtain $\mathbf{u^*}=(\frac{5}{6}, \frac{5}{6}, \frac{5}{6}, \frac{5}{6}, \frac{5}{6})$, which is the greatest optimal solution of problem \eqref{eq13}. Applying Algorithm \ref{alt1}, we have
$$x^1=(\frac{7}{15}, \frac{5}{6}, \frac{5}{6}, \frac{5}{6}, \frac{5}{6}), x^2=(\frac{5}{6}, \frac{5}{6}, \frac{5}{6}, \frac{2}{3}, \frac{5}{6})\in \check{S}(A,b)$$
with $x^1\leq\mathbf{u^*}$ and $x^2\leq\mathbf{u^*}$. Therefore, $x^1$ and $x^2$ are two minimal optimal solutions of problem \eqref{eq13}.

\section{Conclusions}

The contributions of this article include two aspects. One is about minimal solutions of system (\ref{eq2}), in which we first established two necessary and sufficient conditions for a solution of system (\ref{eq2}) to be a minimal one. Then we proved that for every fixed solution of system (\ref{eq2}) there is a minimal solution which is less than or equal to the solution, which means that the solution set of system (\ref{eq2}) is completely determined by its greatest solution and all minimal ones. We also supplied an algorithm for finding minimal solutions of a fixed solution to the system \eqref{eq2} with computational complexity $O(m^3)$ or $O(n^3)$. It should be pointed out that the idea of Algorithm \ref{alt1} originates in \cite{DiNola90}, and some times, Algorithm \ref{alt1} may find more than one minimal solution for a fixed one by changing the permutations of $\{1, 2, \cdots, n\}$ as presented by Example \ref{ex3} even if there is no effective method for finding all minimal solutions of system (\ref{eq2}). The other is to apply minimal solutions of system (\ref{eq2}) to search the optimal solutions of problem \eqref{eq6}. In this part, we related the system (\ref{eq2}) as a constraint to the minimax programming problem, i.e., the problem \eqref{eq6}. We showed that the greatest optimal solution of problem \eqref{eq6} can be obtained by solving a simple single variable optimization problem. Algorithm \ref{al3} served for searching the greatest optimal solution. Algorithm \ref{alt1} was used for finding some minimal optimal solution of problem \eqref{eq6}. In this way, a lot of optimal solutions of problem \eqref{eq6} can be presented with computational complexity $O(m^3)$ or $O(n^3)$ for managing the P2P file sharing system modeled by the system (\ref{eq2}). It is evident that every optimal solution of problem \eqref{eq6} is between a minimal optimal solution and the greatest optimal one since $S^*(A,b)\subseteq S(A,b)$. So that, in the future, we shall develop an effect algorithm for finding all minimal solutions of system (\ref{eq2}) which will further be used for investigating the optimal solutions of problem \eqref{eq6}.


\begin{thebibliography}{9}
\bibitem{Ab06} S. Abbasbandy, E. Babolian and M. Allame, Numerical solution of fuzzy max-min systems, Applied Mathematics and Computation 174 (2006) 1321-1328.
\bibitem{chi} Y. L. Chiu, S. M. Guu, J. J. Yu, Y. K. Wu, A single-variable method for solving min-max programming problem with addition-min fuzzy relational inequalities, Fuzzy Optimization and Decision Making 18 (2019) 433-449.
\bibitem{DiNola90} A. Di Nola, On solving relation equations in Brouwerian lattics, Fuzzy Sets and Systems 34 (1990) 365-376.
\bibitem{DiNola89} A. Di Nola, S. Sessa, W. Pedrycz and E. Sanchez, Fuzzy Relation Equations and Their Applications to Knowledge Engineering, Kluwer Academic Publishers, Dordrecht, Boston/London, 1989.
\bibitem{Guo}F. F. Guo, J. Shen, A smoothing approach for minimizing a linear function subject to fuzzy relation inequalities with addition-min composition, International Journal of Fuzzy Systems 21 (2019) 281-290.
\bibitem{Guu2} S. M. Guu, Y. K. Wu, A linear programming approach for minimizing a linear function subject to fuzzy relational inequalities with addition-min composition, IEEE Transactions on Fuzzy Systems 25 (2017) 985-992.
\bibitem{Guu19} S. M. Guu, Y. K. Wu, Multiple objective optimization for systems with addition-min fuzzy relational inequalities, Fuzzy Optimization and Decision Making 18 (2019) 529-544.
\bibitem{Kl} E. P. Klement, R. Mesiar, E. Pap, Triangular Norms, Kluwer Academic Publishers, Dordrecht, 2000.
\bibitem{Li2012} J. X. Li,  S. J. Yang, Fuzzy relation inequalities about the data transmission mechanism in BitTorrent-like Peer-to-Peer file sharing systems, in: Proceedings of the 2012 9th International Conference on Fuzzy Systems and Knowledge Discover, FSKD 2012, pp. 452-456.
\bibitem{Li2021} M. Li, X. P. Wang, Remarks on minimal solutions of fuzzy relation inequalities with addition-min composition, Fuzzy Sets and Systems 410 (2021) 19-26.
\bibitem{Li2022} M. Li, X. P. Wang, Minimal solutions of fuzzy relation inequalities with addition-min composition and their applications, IEEE Transactions on Fuzzy Systems 31 (2023) 1665-1675.
\bibitem{Lin20} H. T. Lin, X. P. Yang, Dichotomy algorithm for solving weighted min-max programming problem with addition-min fuzzy relation inequalities constraint, Computers and Industrial Engineering 146 (2020) 106537.
\bibitem{LING1965}C. M. Ling, Representation of associative functions, Publ. Math. Debrecen 12 (1965) 189-212.
\bibitem{Martin2017} G. Martin, N. Zuzana, Steady states of max-{\L}ukasiewicz fuzzy systems, Fuzzy Sets and Systems 325 (2017) 58-68.
\bibitem{Peeva04}  K. Peeva, Y. Kyosev, Fuzzy Relational Calculus: Theory, Applications and Software, World Scientific Publishing Company, 2004.
\bibitem{Qiu} J. J. Qiu, X. P.  Yang, Optimization problems subject to addition-{\L}ukasiewicz-product fuzzy relational inequalities with applications in urban sewage treatment systems, Information Sciences 591 (2022) 49-67.
\bibitem{Sanchez76}E. Sanchez, Resolution of composite fuzzy relation equations, Information and Control 30 (1976) 38-48.
\bibitem{So} A. S. Solodovnikov, Systems of linear inequalities, Mir publishers, 1979.
\bibitem{Wu1} Y. K. Wu, Y. L. Chiu, S. M. Guu, Generalized min-max programming problems subject to addition-min fuzzy relational inequalities, Fuzzy Sets and Systems 447 (2022) 22-38.
\bibitem{Wu2} Y. K. Wu, S. M. Guu, An active-set approach to finding a minimal-optimal solution to the min-max programming problem with addition-min fuzzy relational inequalities, Fuzzy Sets and Systems 447 (2022) 39-53.
\bibitem{Wu3} Y. K. Wu, Y. Y. Lur, C. F. Wen, S. J. Lee, Analytical method for solving max-min inverse fuzzy relation, Fuzzy Sets and Systems 440 (2022) 21-41.
\bibitem{Wu} Y. K. Wu, C. F. Wen, Y. T. Hsu, M. X. Wang, Finding minimal solutions to the system of addition-min fuzzy relational inequalities, Fuzzy Optimization and Decision Making 21 (2022) 581-603.
\bibitem{Yang2014} S. J. Yang, An algorithm for minizing a linear objective function subject to the fuzzy relation inequalities with addition-min composition, Fuzzy Sets and Systems 255 (2014) 41-51.
\bibitem{SJYang} S. J. Yang, Some results of fuzzy relation inequalities with addition-min composition, IEEE Transactions on Fuzzy Systems 26 (2018) 239-245.
\bibitem{Yang2017} X. P. Yang, Optimal-vector-based algorithm for solving min-max programming subject to addition-min fuzzy relation inequality, IEEE Transactions on Fuzzy Systems 25 (2017) 1127-1140.
\bibitem{Yang2018} X. P. Yang, H. T. Lin, X. G. Zhou, B. Y. Cao, Addition-min fuzzy relation inequalities with application in BitTorrent-like Peer-to-Peer file sharing system, Fuzzy Sets and Systems 343 (2018) 126-140.
\bibitem{Yang2015} X. P. Yang, X. G. Zhou, B. Y. Cao, Multi-level programming subject to addition-min fuzzy relation inequalities with application in Peer-to-Peer file sharing systems, Journal of Intelligent \& Fuzzy Systems 28 (2015) 2679-2689.
\bibitem{Yang2016} X. P. Yang, X. G. Zhou, B. Y. Cao, Min-max programming problem subject to addition-min fuzzy relation inequalities, IEEE Transactions on Fuzzy Systems 24 (2016) 111-119.
\end{thebibliography}
\end{document}